\documentclass[11pt,a4paper]{article}
\pdfoutput=1
\usepackage{lmodern}
\usepackage[normalem]{ulem}
\usepackage{titlesec}
\titleformat{\section}
  {\normalfont\fontsize{16pt}{10pt}\selectfont\scshape}
  {\thesection}
  {0.6em}
  {}
\titleformat{\subsection}
  {\normalfont\fontsize{14pt}{8pt}\selectfont\scshape}
  {\thesubsection}
  {0.6em}
  {}

\usepackage{amsmath,amssymb,amsfonts,amsthm,mathtools}
\usepackage{mathrsfs}
\usepackage{color}
\usepackage[margin=2.435cm]{geometry}
\usepackage{bbm}
\usepackage{graphicx}
\usepackage{enumitem}
\usepackage{hyperref}

\definecolor{darkred}{rgb}{0.8,0.1,0.1}
\hypersetup{
	colorlinks=true,         
	linkcolor=darkred,
	citecolor=blue,
}

\newtheoremstyle{scstyle}%
  {0.4em}   
  {0.4em}   
  {}      
  {}      
  {\fontsize{12pt}{8pt}\selectfont\scshape} 
  {.}     
  {.5em}  
  {}      

\theoremstyle{scstyle}
\newtheorem{theo}{Theorem}[section]
\newtheorem{lem}[theo]{Lemma}
\newtheorem{propo}[theo]{Proposition}

\newtheorem{cor}[theo]{Corollary}

\theoremstyle{scstyle}
\newtheorem{defi}[theo]{Definition}

\newtheorem{exx}[theo]{Example}
\newenvironment{ex}
{\pushQED{\qed}\exx}
{\popQED\endexx}

\newtheorem{remm}[theo]{Remark}
\newenvironment{rem}
{\pushQED{\qed}\remm}
{\popQED\endremm}

\numberwithin{equation}{section}

\def\nn{\nonumber}
\def\A{\mathcal{A}}
\def\Li{\mathrm{Li}}
\def\I{\mathcal{I}}
\def\bbC{\mathbb{C}}
\def\bbN{\mathbb{N}}
\def\ad{\mathrm{ad}}

\def\sk{\vspace{2mm}}

\title{Relations between multiple zeta values and delta values from Drinfeld's associator series}

\author{%
Cameron James Deverall Kemp\vspace{4mm}\\
{\small School of Mathematical Sciences, University of Nottingham,}\\
{\small University Park, Nottingham NG7 2RD, United Kingdom.}\vspace{4mm}\\
{\small \begin{tabular}{ll}
Email: & \href{mailto:cameron.kemp@nottingham.ac.uk}{\texttt{cameronkemp31@gmail.com}}
\vspace{2mm}
\end{tabular}
}
}

\begin{document}
\maketitle

\begin{abstract}
\noindent It is shown that novel relations between multiple zeta values and single-variable multiple polylogarithms at 1/2 (delta values) can be derived by comparing two distinct, yet a priori equal, series formulae for the Drinfeld associator (from the Knizhnik-Zamolodchikov connection). In particular, we demonstrate that two new relations are found by comparing the fifth-order terms of each series formula. 
\end{abstract}



\tableofcontents

\section{Introduction and summary}
In \cite{Drinfeld90}, Drinfeld introduced a $\bbC[[\hbar]]$-coefficient power series $\Phi_\mathrm{KZ}(A,B)$ which, given a complex Lie algebra $\mathfrak{g}$ and a quadratic Casimir $t\in(\mathrm{Sym}^2\mathfrak{g})^\mathfrak{g}$, deforms the cocommutative Hopf algebra $U\mathfrak{g}$ into a quasitriangular quasi-Hopf algebra $U\mathfrak{g}[[\hbar]]$ (a Drinfeld-Jimbo quantum group) where the quasitriangular structure is given by $R=e^{i\pi\hbar t}$ and the associator is given by $\Phi=\Phi_\mathrm{KZ}(t_{12},t_{23})$. The usual method of constructing an explicit series formula for $\Phi_\mathrm{KZ}(A,B)$ is presented in \cite[Section XIX.6]{Kassel}, see \cite[Theorem A.9]{LeMu} for the explicit series formula itself. In particular, this method leads to the coefficients being MZVs, see Definition \ref{def:MZV}. MZVs are special cases of multiple polylogarithms, see Definition \ref{def:multiple polylogarithm}.  
\sk

Multiple polylogarithms appear naturally in a wide range of topics, e.g.: knot theory (where they appear as MZVs) \cite{LeMu,phi4,counterterms}, algebraic and hyperbolic geometry \cite{Goncharov,Goncharov1,Goncharov2}, the evaluation of multi-loop Feynman diagrams in perturbative quantum field theory \cite{Bogner,Broadhurst,Heller} and quantum statistical thermodynamics (as polylogarithms) \cite{Guan,Howard}. Throughout their study, relations between multiple polylogarithms play a central role \cite{LeMu0,SCstar,Polylog,SVoMP}. 
\sk

Recently, an alternative method of constructing an explicit series formula for $\Phi_\mathrm{KZ}(A,B)$ was presented in \cite[Lemma 19 and Theorem 20]{BRW}. Bordemann, Rivezzi and Weigel (BRW) state in the introduction that the objective of their work was to provide a pedagogical approach to the work of Drinfeld \cite{Drinfeld90} and not contain any new result. However, we will show in Proposition \ref{propo:iterated integral is sum of deltas} that their method provides a series formula whose coefficients evaluate to delta values, another special case of a multiple polylogarithm, see Definition \ref{def:single-variable multiple polylogarithm}. In other words, the new result implicit in their work was an infinitude of nontrivial relations between MZVs and delta values. The present paper finds two novel such relations at order $\hbar^5$ by making use of iterated integral and parallel transport techniques; the relevance being that, as mentioned above, relations between multiple polylogarithms are central to their study.
\sk

For instance, equating the two series formulae reveals, at order $\hbar^2$, the relation 
\begin{equation}\label{eq:Li_2(1/2)}
\zeta_2=2\delta_2+(\ln2)^2\quad.
\end{equation}
Of course, the relation \eqref{eq:Li_2(1/2)} is commonly known; in fact, it is \textbf{Euler's dilogarithm formula}, see \cite[(1.7)]{Lewin}. The novelty of our arguments is that we can (re-)derive relations such as \eqref{eq:Li_2(1/2)} independently from the standard route via clever manipulation of integral expressions, nested sums, etc. As another example, at order $\hbar^3$, we find
\begin{equation}\label{eq:delta(2,1)}
\zeta_3=\,\delta_{2,1}+\tfrac{1}{2}(\ln2)^3+\delta_2\ln2+\delta_3\quad.
\end{equation}
The relation \eqref{eq:delta(2,1)} is far less well known and the only work we could find it in was \cite[(7.3)]{SVoMP}, see \eqref{eq:Holder convolution} with $n=1$. At order $\hbar^4$, we have the following two relations:
\begin{subequations}
\begin{alignat}{2}
\zeta_4&=\delta_{2,1,1}+\tfrac{1}{6}(\ln2)^4+\tfrac{1}{2}\delta_2(\ln2)^2+\delta_3\ln2+\delta_4\quad&&,\label{eq:Holder with n=2}\\
\tfrac{1}{4}\zeta_4&=2\delta_{3,1}+\delta_{2,1}\ln2+\tfrac{1}{4}(\ln2)^4\quad&&.\label{eq:crazy relation}
\end{alignat}
\end{subequations}
The relation \eqref{eq:Holder with n=2} is simply \eqref{eq:Holder convolution} with $n=2$. The other relation \eqref{eq:crazy relation} is a rich interplay of various relations from \cite{SVoMP} together with \cite[(62)]{BBB96} and \cite[(1.14)]{Lewin}. The idea is that the higher one goes in the order of $\hbar$ the more relations one (re-)discovers and the more esoteric those relations become. In particular, at order $\hbar^5$, we find:
\begin{theo}
\begin{subequations}
\begin{alignat}{3}
\zeta_5=&\,\delta_{2,1,1,1}+\tfrac{1}{24}(\ln2)^5+\tfrac{1}{6}\delta_2(\ln2)^3+\tfrac{1}{2}\delta_3(\ln2)^2+\delta_4\ln2+\delta_5\quad&&,\label{eq:Holder n=3}\\
\zeta_{4,1}=&\,\delta_{4,1}+\delta_{3,1}\ln2+\tfrac{1}{2}\delta_{2,1}(\ln2)^2+\delta_{3,1,1}+\delta_{2,1,1}\ln2+\tfrac{1}{12}(\ln2)^5\quad&&,\label{eq:delta_4,1 relation}\\
3\zeta_{4,1}+\zeta_{3,2}=&\,\delta_{3,2}+3\delta_{4,1}+\left(\tfrac{1}{2}\delta_2^2+\delta_{3,1}\right)\ln2+\left[\tfrac{1}{2}(\ln2)^2+\zeta_2\right]\delta_{2,1}-3\delta_{3,1,1}\nn\\&-2\delta_{2,2,1}-\delta_{2,1,2}+\tfrac{1}{4}\zeta_2(\ln2)^3\quad&&.\label{eq:delta_3,2 relation}
\end{alignat}
\end{subequations}
\end{theo}
The relation \eqref{eq:Holder n=3} is yet another instance of \eqref{eq:Holder convolution} ($n=3$) whereas the relations \eqref{eq:delta_4,1 relation} and \eqref{eq:delta_3,2 relation} are, to the best of our knowledge, completely new and never-before-seen, see Remark \ref{rem:genuinely new}.

\subsection*{Plan of the paper}
Section \ref{sec:preliminaries} recalls some basic definitions and integral expressions as well as various relations that will be used throughout. In particular, Subsection \ref{subsec:Multiple polylogs} recalls the definition of: the Riemann zeta function, the gamma function, (single-variable multiple) polylogarithms, as well as their various integral expressions and relevant special values. Subsection \ref{subsec:iterated integrals} recalls the notion of Chen's iterated integrals and the Kontsevich iterated integral expression of single-variable multiple polylogarithms as in \cite[directly below (4.13)]{SVoMP} together with the inversion property and shuffle product property that iterated integrals satisfy; the MZV duality and shuffle relations will be derived from these properties, respectively. 
\sk

Section \ref{sec:associator from PT} recalls the two distinct methods of constructing the explicit power series formula (in $\hbar$) for the Drinfeld associator $\Phi_\mathrm{KZ}(A,B)$ as the limit $\lim_{\varepsilon\to0}$ of the non-singular factor of a parallel transport with respect to the formal Knizhnik-Zamolodchikov connection from $0<\varepsilon<\tfrac{1}{4}$ to $\frac{3}{4}<1-\varepsilon<1$. In Subsection \ref{subs:MZV PT}, we recall the first method as found in \cite[Section XIX.6]{Kassel}. The explicit result of the first method yields \cite[Theorem A.9]{LeMu} and we restate this ``MZV-formula" in Corollary \ref{cor:Phi with MZV coeffs}. In Subsection \ref{subsec:BRW PT}, we recall the second method as found in \cite[Lemma 19 and Theorem 20]{BRW}. The explicit result of the second method is not known in the literature but we provide this ``delta-formula" in Proposition \ref{propo:iterated integral is sum of deltas}.
\sk

Section \ref{sec:Fifth-order expansion} carries out the expansion of each series formula up to the fifth-order term. It will be noticed that the fifth-order expansion of the MZV-formula \eqref{eq:fifth order MZV expansion} is much simpler than that of the delta-formula in Subsection \ref{subsec:DV expansion}. Subsection \ref{subsec:relations} expresses, order-by-order, the various relations imposed by the fact that these two series formulae are equal. 

\section{Preliminaries}\label{sec:preliminaries}
Many of the definitions and expressions appearing in this section are taken from \cite[Appendix XIX.11]{Kassel}, \cite{Polylog} and \cite{SVoMP}.
\subsection{Multiple polylogarithms}\label{subsec:Multiple polylogs}
See \cite[(1.3)]{SVoMP} for the following definition.
\begin{defi}
The \textbf{Riemann zeta function} is defined as 
\begin{equation}\label{eq:Riemann zeta as infinite sum}
\zeta_s:=\sum_{n=1}^\infty\frac{1}{n^s}
\end{equation}
for $\mathrm{Re}(s)>1$ and analytically continued elsewhere.
\end{defi}
The Riemann zeta function admits an integral expression (for $\mathrm{Re}(s)>1$) and for this we recall the \textbf{gamma function} which is defined as
\begin{equation}\label{eq:Gamma function}
\Gamma(z):=\int_0^\infty t^{z-1}e^{-t}dt
\end{equation}
for $\mathrm{Re}(z)>0$ and analytically continued using the identity $\Gamma(z)=\frac{\Gamma(z+1)}{z}$. The gamma function is thus a meromorphic function with simple poles at $z\in\mathbb{Z}^{\leq0}$. The analytic continuation of the Riemann zeta function uses \textbf{Riemann's reflection formula} $\zeta_s=2^s\pi^{s-1}\sin(\frac{\pi s}{2})\Gamma(1-s)\zeta_{1-s}$ hence the Riemann zeta function is a meromorphic function with a simple pole at $s=1$.
\begin{rem}
The gamma function is said to be an extension of the factorial function to the complex numbers because, for $k\in\bbN$, 
\begin{equation}
\Gamma(k+1):=\int_0^\infty t^ke^{-t}dt=k!\quad.
\end{equation}
\end{rem}
Now we can derive the integral representation of the Riemann zeta function as\footnote{Recall that the gamma function possesses no zeroes (which can be shown using the identity $\Gamma(z+1)=z\Gamma(z)$ and \textbf{Euler's reflection formula} $\Gamma(1-z)\Gamma(z)=\frac{\pi}{\sin(\pi z)}$) thus its reciprocal $\frac{1}{\Gamma(s)}$ is an entire function.} 
\begin{equation}\label{eq:zeta as integral}
\frac{1}{\Gamma(s)}\int_0^\infty\frac{x^{s-1}}{e^x-1}dx=\frac{1}{\Gamma(s)}\int_0^\infty x^{s-1}\sum_{n=1}^\infty e^{-nx}dx=\sum_{n=1}^\infty\frac{1}{n^s}\frac{1}{\Gamma(s)}\int_0^\infty t^{s-1}e^{-t}dt=\zeta_s
\end{equation}
where the first equality uses the Taylor expansion $\frac{1}{e^x-1}=e^{-x}(1-e^{-x})^{-1}=\sum_{n=1}^\infty e^{-nx}$, the second equality uses the substitution $nx=t$ and the last equality uses the definition \eqref{eq:Gamma function} of the gamma function with argument $z=s$. 

Let us recall the following extension of the Riemann zeta function \eqref{eq:Riemann zeta as infinite sum}.
\begin{defi}
The \textbf{polylogarithm} is defined as
\begin{equation}\label{eq:definition of polylogarithm}
\Li_s(z):=\sum_{n=1}^\infty\frac{z^n}{n^s}
\end{equation}
for $(s\in\bbC,|z|<1)$ and analytically continued elsewhere. We refer to $s\in\bbC$ as the \textbf{order}.
\end{defi}
\begin{rem}
The polylogarithm \eqref{eq:definition of polylogarithm} is indeed an extension of the Riemann zeta function \eqref{eq:Riemann zeta as infinite sum} because $\Li_s(1)=\zeta_s$ but it is \textit{also} an extension of the natural logarithm function,
\begin{equation}\label{eq:Li_1(z)=-ln(1-z)}
\Li_1(z):=\sum_{n=1}^\infty\frac{z^n}{n}=-\ln(1-z)\quad.
\end{equation}
\end{rem}
\begin{ex}\label{ex:polylog of -1}
Using \eqref{eq:Li_1(z)=-ln(1-z)} we have
\begin{subequations}
\begin{equation}\label{eq:Li_1(-1)}
\Li_1(-1)=-\ln2\quad,
\end{equation}
whereas for $\mathrm{Re}(s)>1$  we have \cite[Example 9.2]{SVoMP}
\begin{equation}\label{eq:Li_s(-1)}
\Li_s(-1)=(2^{1-s}-1)\zeta_s\quad.
\end{equation}
\end{subequations}
\end{ex}
The polylogarithm admits the following integral expressions for $(s\in\bbC,|z|<1)$:
\begin{enumerate}
\item[(i)] We have
\begin{subequations}
\begin{equation}
\Li_s(z)=\int_0^z\sum_{n=1}^\infty\frac{r^{n-1}}{n^{s-1}}\mathrm{d}r=\int_0^z\frac{\Li_{s-1}(r)}{r}\mathrm{d}r\quad.
\end{equation}
\item[(ii)] The \textbf{Bose-Einstein integral}\footnote{Derived analogously to \eqref{eq:zeta as integral}.},
\begin{equation}\label{eq:Bose-Einstein integral}
\Li_s(z)=\frac{z}{\Gamma(s)}\int_0^\infty\frac{r^{s-1}}{e^r-z}\mathrm{d}r\quad.
\end{equation}
\item[(iii)] Substituting $r=-\ln u$ into \eqref{eq:Bose-Einstein integral} gives
\begin{equation}\label{eq:POLYlog integral}
\Li_s(z)=\frac{z}{\Gamma(s)}\int_0^1\frac{\left(-\ln u\right)^{s-1}}{1-uz}\mathrm{d}u\quad.
\end{equation}
\end{subequations}
\end{enumerate}
In addition to the polylogarithm, there exists another generalisation of the Riemann zeta function.
\begin{defi}\label{def:MZV}
Given $k\in\bbN^*$, the \textbf{multiple zeta function} is defined as 
\begin{equation}\label{eq:multiple zeta function}
\zeta_{s_1,\ldots,s_k}:=\sum_{n_1>n_2>\cdots>n_k\geq1}\frac{1}{n_1^{s_1}\cdots n_k^{s_k}}
\end{equation}
for $\{\sum_{j=1}^i\mathrm{Re}(s_j)>i\,|\,1\leq i\leq k\}$ and analytically continued elsewhere. If $s_1\in\bbN^{\geq2}$ and $s_2,\ldots,s_k\in\bbN^*$ then we call \eqref{eq:multiple zeta function} a \textbf{multiple zeta value} (MZV). 
\end{defi}
We will need the ``self-dual evaluation" \cite[Example 6.3 with $n=1$]{SVoMP}
\begin{subequations}
\begin{equation}\label{eq:zeta_3,1}
\zeta_{3,1}=\frac{1}{4}\zeta_4
\end{equation}
as well as \cite[Example 2.1 with $m=4$]{SVoMP}
\begin{equation}\label{eq:zeta_4,1}
\zeta_{4,1}=2\zeta_5-\zeta_2\zeta_3
\end{equation}
and the ``stuffle identity" \cite[Top of subsection 3.1 with $(s=2,t=3)$]{SVoMP}
\begin{equation}\label{eq:stuffle identity}
\zeta_2\zeta_3=\zeta_{2,3}+\zeta_{3,2}+\zeta_5\quad.
\end{equation} 
\end{subequations}
See \cite[(1.4)]{SVoMP} for the following generalisation of both the polylogarithm \eqref{eq:definition of polylogarithm} and the multiple zeta function \eqref{eq:multiple zeta function}.
\begin{defi}\label{def:multiple polylogarithm}
Given $k\in\bbN^*$, the \textbf{multiple polylogarithm} is defined as 
\begin{equation}\label{eq:multiple polylogarithm}
\Li_{s_1,\ldots,s_k}(z_1,\ldots,z_k):=\sum_{n_1>n_2>\cdots>n_k\geq1}\frac{z_1^{n_1}\cdots z_k^{n_k}}{n_1^{s_1}\cdots n_k^{s_k}}
\end{equation}
for $\{(s_i\in\bbC,|z_i|<1)\,|\,1\leq i\leq k\}$ and analytically continued elsewhere.
\end{defi}
\begin{ex}
For $n\in\bbN^*$, we set $\{z\}^n:=(\overbrace{z,\ldots,z}^{n})$. Generalising \eqref{eq:Li_1(-1)}, \cite[(62)]{BBB96} gives
\begin{equation}\label{eq:Li_1,...,1(-1,...,-1)}
\Li_{\{1\}^n}(\{-1\}^n)=(-1)^n\sum_{\sum_{k=1}^nkj_k=n}\,\prod_{k=1}^n\frac{1}{j_k!}\left(\frac{-\Li_k\big((-1)^k\big)}{k}\right)^{j_k}\quad.
\end{equation}
In particular, we will need:
\begin{subequations}\label{subeq:Li_1^4,5(-1)}
\begin{equation}\label{eq:Li_1^4(-1)}
\Li_{\{1\}^4}(\{-1\}^4)=\frac{(\ln2)^4}{24}-\frac{\zeta_2(\ln2)^2}{4}+\frac{\zeta_2^2}{8}+\frac{\zeta_3\ln2}{4}-\frac{\zeta_4}{4}
\end{equation}
and
\begin{equation}\label{eq:Li_1^5(-1)}
\Li_{\{1\}^5}(\{-1\}^5)=-\frac{3\zeta_5}{16}+\frac{\zeta_4\ln2}{4}-\frac{\zeta_3(\ln2)^2}{8}+\frac{\zeta_2\zeta_3}{8}+\frac{\zeta_2(\ln2)^3}{12}-\frac{\zeta_2^2\ln2}{8}-\frac{(\ln2)^5}{120}\quad,
\end{equation}
\end{subequations}
where we used Example \ref{ex:polylog of -1} for both calculations.
\end{ex}

\begin{defi}\label{def:single-variable multiple polylogarithm}
Setting $z_1=z$ and $z_2=\ldots=z_k=1$, we refer to 
\begin{equation}\label{eq:single-variable multiple polylogarithm}
\Li_{s_1,\ldots,s_k}(z):=\Li_{s_1,\ldots,s_k}(z,\{1\}^{k-1}):=\sum_{n_1>n_2>\cdots>n_k\geq1}\frac{z^{n_1}}{n_1^{s_1}\cdots n_k^{s_k}}
\end{equation}
as a \textbf{single-variable multiple polylogarithm}. We refer to: $k\in\bbN^*$ as the \textbf{depth}, $s_i$ as the \textbf{$i^\mathrm{th}$-order} and $\sum_{i=1}^ks_i$ as the \textbf{weight}. We are particularly interested in the case $z=\tfrac{1}{2}$ and $s_1,\ldots,s_k\in\bbN^*$ thus we adopt the notation of \cite[(2.2)]{SVoMP} and define the \textbf{delta value},
\begin{equation}\label{eq:definition of delta}
\delta_{s_1,\ldots,s_k}:=\Li_{s_1,\ldots,s_k}\left(\tfrac{1}{2}\right)\quad.
\end{equation}
\end{defi}
\begin{rem}
The single-variable multiple polylogarithm \eqref{eq:single-variable multiple polylogarithm} is clearly a generalisation of the polylogarithm \eqref{eq:definition of polylogarithm} and an extension of the multiple zeta function \eqref{eq:multiple zeta function}. Furthermore, something analogous to \eqref{eq:Li_1(z)=-ln(1-z)} holds, i.e.
\begin{equation}\label{Li_(1,...,1)(z)}
\Li_{\{1\}^k}(z)=\frac{\big(-\ln(1-z)\big)^k}{k!}\quad,
\end{equation}
in particular, 
\begin{equation}\label{eq:delta_(1^k)}
\delta_{\{1\}^k}=\frac{(\ln2)^k}{k!}\quad.
\end{equation}
\end{rem}
As in \cite[(11.18) and (11.19) of Appendix XIX.11]{Kassel}, we mention the following integral expressions:
\begin{equation}\label{eq:(11.18/19) of Kassel}
\Li_{s_1+1,s_2\ldots,s_k}(z)=\int_0^z\frac{\Li_{s_1,s_2\ldots,s_k}(r)}{r}dr\qquad,\qquad\Li_{1,s_1,\ldots,s_k}(z)=\int_0^z\frac{\Li_{s_1,\ldots,s_k}(r)}{1-r}dr\quad.
\end{equation}
\subsection{Some known
relations between MZVs and delta values}
This short subsection collects a few known relations between MZVs and delta values; most of the relations in Subsection \ref{subsec:relations} will be shown to be derivative from the relations of this subsection whereas the last two relations in Subsection \ref{subsec:relations} will be shown to be independent. We recall \textbf{Landen's trilogarithm formula} from \cite[(1.14)]{Lewin},
\begin{equation}\label{eq:Landen's trilog formula}
\delta_3=\frac{7\zeta_3}{8}-\frac{\zeta_2\ln2}{2}+\frac{(\ln2)^3}{6}\quad.
\end{equation}
To the author’s
knowledge, formulae for delta values purely in terms of zeta values and $\ln2$ are not known
for higher orders. That being said, we also have 
\begin{equation}
\delta_4=\frac{\zeta_4}{4}+\frac{\zeta_2(\ln2)^2}{4}-\frac{(\ln2)^4}{24}-\frac{1}{2}\Li_{3,1}(-1,-1)
\end{equation}
as in \cite[page 291]{BBG}. The following lemma recalls the notion of \textbf{Hölder convolution} from \cite[(7.3)]{SVoMP}.
\begin{lem}
For $n\in\bbN$,
\begin{equation}\label{eq:Holder convolution}
\zeta_{n+2}=\delta_{2,\{1\}^n}+\sum_{r=1}^{n+2}\frac{(\ln2)^{n+2-r}}{(n+2-r)!}\delta_r\quad.
\end{equation}
\end{lem}
\begin{proof}
Substitute $k=b_1=1$, $s_1=n+2$ and $p=q=2$ into the formula \cite[(7.1)]{SVoMP} to get $\zeta_{n+2}=\delta_{2,\{1\}^n}+\sum_{r=1}^{n+2}\delta_{\{1\}^{n+2-r}}\delta_r$ then simply recall \eqref{eq:delta_(1^k)}.
\end{proof}
Notice that Euler's dilogarithm formula \eqref{eq:Li_2(1/2)} can be seen as the case $n=0$ of \eqref{eq:Holder convolution}. Substituting Landen's trilogarithm formula \eqref{eq:Landen's trilog formula} and Euler's dilogarithm formula \eqref{eq:Li_2(1/2)} into \eqref{eq:delta(2,1)} gives
\begin{equation}\label{delta_2,1}
\delta_{2,1}=\frac{\zeta_3}{8}-\frac{(\ln2)^3}{6}\quad.
\end{equation}
Using both \cite[(2.3) with (6.8)]{SVoMP} and \eqref{subeq:Li_1^4,5(-1)} gives: 
\begin{subequations}\label{subeq:delta_2,2 and delta_1,2,2}
\begin{alignat}{2}
\delta_{2,2}~=&~\frac{(\ln2)^4}{24}-\frac{\zeta_2(\ln2)^2}{4}+\frac{\zeta_2^2}{8}+\frac{\zeta_3\ln2}{4}-\frac{\zeta_4}{4}&&,\label{eq:delta_2,2}\\
\delta_{1,2,2}~=&\frac{3\zeta_5}{16}-\frac{\zeta_4\ln2}{4}+\frac{\zeta_3(\ln2)^2}{8}-\frac{\zeta_2\zeta_3}{8}-\frac{\zeta_2(\ln2)^3}{12}+\frac{\zeta_2^2\ln2}{8}+\frac{(\ln2)^5}{120}\quad&&.\label{eq:delta_1,2,2}
\end{alignat}
\end{subequations}

\subsection{Iterated integrals}\label{subsec:iterated integrals}
Following \cite[Appendix XIX.11]{Kassel}, suppose we have $\bbC$-valued 1-forms $\{\omega_i:=f_i(s)ds\}_{1\leq i\leq n}$ defined over the real interval $[a,b]\subset\mathbb{R}\,$. Define the iterated integral $\int_a^b\omega_1\cdots\omega_n$ inductively by
\begin{subequations}
\begin{equation}
\int_a^b\omega_1=\int_a^bf_1(s)ds
\end{equation}
and
\begin{equation}
\int_a^b\omega_1\cdots\omega_n=\int_a^bf_1(s)\left(\int_a^s\omega_2\cdots\omega_n\right)ds\quad.
\end{equation}
\end{subequations}
It is straightforward to prove that iterated integrals satisfy the \textbf{inversion property}
\begin{subequations}
\begin{equation}\label{eq:inversion for iterated integrals}
\int_a^b\omega_1\cdots\omega_n=(-1)^n\int^a_b\omega_n\cdots\omega_1\quad.
\end{equation}
For $m,n\in\bbN^*$, we say a permutation $\sigma\in\mathrm{S}_{m+n}$ is an \textbf{$(m,n)$-shuffle} of $(1,\ldots,m,m+1,\ldots,m+n)$ if $\sigma(1)<\cdots<\sigma(m)$ and $\sigma(m+1)<\cdots<\sigma(m+n)$. We denote the set of $(m,n)$-shuffles as $\mathrm{S}_{m,n}$. It is also straightforward to prove that iterated integrals satisfy the \textbf{shuffle product property}
\begin{equation}\label{eq:shuffle prod for iterated integrals}
\int_a^b\omega_1\cdots\omega_m\int_a^b\omega_{m+1}\cdots\omega_{m+n}=\sum_{\sigma\in\mathrm{S}_{m,n}}\int_a^b\omega_{\sigma(1)}\cdots\omega_{\sigma(m+n)}\quad.
\end{equation}
\end{subequations}
\begin{rem}
Define the particular 1-forms:
\begin{equation}\label{eq:particular 1-forms}
\Omega_0:=\frac{ds}{s}\qquad,\qquad\Omega_1:=\frac{ds}{s-1}\quad.
\end{equation}
Then, for $g\in\bbN^*$ and $p_1,q_1,\ldots,p_g,q_g\in\bbN^*$, we can use \eqref{eq:(11.18/19) of Kassel} to notice that
\begin{equation}\label{eq:iterated integral expression for multiple zeta}
\int_0^1\Omega_0^{p_1}\Omega_1^{q_1}\cdots\Omega_0^{p_g}\Omega_1^{q_g}=(-1)^{\sum_{i=1}^gq_i}\zeta_{p_1+1,\{1\}^{q_1-1},\ldots,p_g+1,\{1\}^{q_g-1}}\quad.
\end{equation}
In particular, we can use the inversion property \eqref{eq:inversion for iterated integrals} and integration by substitution,
\begin{align}
\int_0^1\Omega_0^{p_1}\Omega_1^{q_1}\cdots\Omega_0^{p_g}\Omega_1^{q_g}&\,=(-1)^{\sum_{i=1}^g(p_i+q_i)}\int_1^0\Omega_1^{q_g}\Omega_0^{p_g}\cdots\Omega_1^{q_1}\Omega_0^{p_1}\nn\\
&\,=(-1)^{\sum_{i=1}^g(p_i+q_i)}\int_0^1\Omega_0^{q_g}\Omega_1^{p_g}\cdots\Omega_0^{q_1}\Omega_1^{p_1}\quad,
\end{align}
to derive the \textbf{MZV duality}, 
\begin{equation}\label{eq:MZV duality}
\zeta_{p_1+1,\{1\}^{q_1-1},\ldots,p_g+1,\{1\}^{q_g-1}}=\,\zeta_{q_g+1,\{1\}^{p_g-1},\ldots,q_1+1,\{1\}^{p_1-1}}\quad.
\end{equation}
\end{rem}
\begin{ex}
We mention the following relevant instances of MZV duality \eqref{eq:MZV duality}:
\begin{equation}\label{eq:relevant MZV dualities}
\zeta_3=\zeta_{2,1}~,\quad\zeta_4=\zeta_{2,1,1}~,\quad\zeta_5=\zeta_{2,1,1,1}~,\quad\zeta_{4,1}=\zeta_{3,1,1}~,\quad\zeta_{3,2}=\zeta_{2,2,1}~,\quad\zeta_{2,3}=\zeta_{2,1,2}~.
\end{equation}
\end{ex}
More generally, for $k=r\in\bbN^*$, let us substitute $b_1=\ldots=b_r=\frac{1}{z}$ into the formula \cite[below (4.13)]{SVoMP} to get 
\begin{equation}\label{eq:single variable mult polylog as iterated int}
\Li_{s_1,\ldots,s_r}(z)=\prod_{j=1}^r\int_0^{y_{j-1}}\frac{\left(\ln\frac{y_{j-1}}{y_j}\right)^{s_j-1}}{\Gamma(s_j)(\tfrac{1}{z}-y_j)}dy_j\quad.
\end{equation}
We will need this ``Kontsevich iterated integral formula" for the single-variable multiple polylogarithm in Subsection \ref{subsec:BRW PT}. 

\section{Two constructions of Drinfeld's KZ associator}\label{sec:associator from PT}
Following \cite[Subsection 2.2]{BRW}, given a unital associative $\bbC$-algebra $\A$ and two elements $A,B\in\A$, we define the (flat) formal \textbf{Knizhnik-Zamolodchikov connection} 
\begin{equation}
\Gamma(A,B):=\hbar\left(\frac{A}{x}+\frac{B}{x-1}\right)dx
\end{equation}
on $(0,1)\subset\mathbb{R}$. The parallel transport from $0<\varepsilon<\tfrac{1}{4}$ to $\frac{3}{4}<1-\varepsilon<1$ along the affine path $c_\varepsilon(s):=\varepsilon+s(1-2\varepsilon)$ is given by
\begin{equation}\label{eq:c_epsilon PT}
W^{c_\varepsilon}=1+\sum_{r=1}^\infty\hbar^r\sum_{i_1,\ldots,i_r=0}^1\int_\varepsilon^{1-\varepsilon}\frac{\mathrm{d}u_1}{u_1-i_1}\cdots\frac{\mathrm{d}u_r}{u_r-i_r}\A_{i_1}\cdots\A_{i_r}
\end{equation}
where $\A_0:=A$ and $\A_1:=B$. 
\begin{rem}\label{rem:at least diverging terms}
If $[A,B]=0$ we can simply exponentiate the integration of the connection,
\begin{equation}\label{eq:troublesome terms exposed}
W^{c_\varepsilon}=\exp\left[\hbar\int_\varepsilon^{1-\varepsilon}\left(\frac{A}{x}+\frac{B}{x-1}\right)dx\right]=e^{\hbar\ln(\varepsilon)B}e^{\hbar\ln(1-\varepsilon)(A-B)}e^{-\hbar\ln(\varepsilon)A}
\end{equation}
and this clearly diverges as $\varepsilon\to0$.
\end{rem}
Remark \ref{rem:at least diverging terms}, particularly \eqref{eq:troublesome terms exposed}, makes it clear that we need to extract \textit{at least} the terms $e^{\hbar\ln(\varepsilon)B}$ and $e^{-\hbar\ln(\varepsilon)A}$ in order to acquire something that is \textit{potentially} finite. In other words, we want to show that 
\begin{equation}\label{eq:potentially finite Drinfeld}
\Phi(A,B):=\lim_{\varepsilon\to0}\Phi_\varepsilon(A,B):=\lim_{\varepsilon\to0}\left(e^{-\hbar\ln(\varepsilon)B}W^{c_\varepsilon}e^{\hbar\ln(\varepsilon)A}\right)
\end{equation}
is actually finite.

\subsection{With coefficients as MZVs}\label{subs:MZV PT}
We use \eqref{eq:particular 1-forms} to rewrite \eqref{eq:c_epsilon PT} as 
\begin{align}
W^{c_\varepsilon}=1&+\sum_{r=1}^\infty\hbar^r\sum_{i_1,\ldots,i_r=0}^1\int_\varepsilon^{1-\varepsilon}\Omega_{i_1}\cdots\Omega_{i_r}\A_{i_1}\cdots\A_{i_r}\nn\\
=1&+\hbar\int_\varepsilon^{1-\varepsilon}(\Omega_0A+\Omega_1B)+\sum_{r=2}^\infty\hbar^r\sum_{i_2,\ldots,i_r=0}^1\int_\varepsilon^{1-\varepsilon}\Omega_1\Omega_{i_2}\cdots\Omega_{i_r}B\A_{i_2}\cdots\A_{i_r}\nn\\&+\sum_{r=2}^\infty\hbar^r\sum_{i_2,\ldots,i_{r-1}=0}^1\int_\varepsilon^{1-\varepsilon}\Omega_0\Omega_{i_2}\cdots\Omega_{i_{r-1}}\Omega_0A\A_{i_2}\cdots\A_{i_{r-1}}A\nn\\
&+\sum_{r=2}^\infty\hbar^r\sum_{i_2,\ldots,i_{r-1}=0}^1\int_\varepsilon^{1-\varepsilon}\Omega_0\Omega_{i_2}\cdots\Omega_{i_{r-1}}\Omega_1A\A_{i_2}\cdots\A_{i_{r-1}}B\quad.\label{eq:c_epsilon compact PT}
\end{align}
It is obvious that $\lim_{\varepsilon\to0}\int_\varepsilon^{1-\varepsilon}\Omega_{i_1}\cdots\Omega_{i_r}$ is divergent if $i_r=0$ or, by inversion \eqref{eq:inversion for iterated integrals}, if $i_1=1$. Otherwise, we have 
\begin{align}
\sum_{r=2}^\infty\hbar^r\sum_{i_2,\ldots,i_{r-1}=0}^1\int_0^1\Omega_0\Omega_{i_2}\cdots\Omega_{i_{r-1}}\Omega_1=&\sum _{\begin{smallmatrix}g\geq1\\p_1,\ldots,p_g\geq1\\q_1,\ldots,q_g\geq1\end{smallmatrix}}\hbar^{\sum_{i=1}^g(p_i+q_i)}\int_0^1\Omega_0^{p_1}\Omega_1^{q_1}\cdots\Omega_0^{p_g}\Omega_1^{q_g}\label{eq:convergent iterated integral}\\
=&\sum _{\begin{smallmatrix}g\geq1\\p_1,..,p_g\geq1\\q_1,..,q_g\geq1\end{smallmatrix}}\hbar^{\sum_{i=1}^gp_i+q_i}(-1)^{\sum_{i=1}^gq_i}\zeta_{p_1+1,\{1\}^{q_1-1},\ldots,p_g+1,\{1\}^{q_g-1}}\nn
\end{align}
where the second equality used \eqref{eq:iterated integral expression for multiple zeta}. Following the discussion after \cite[(6.9) of Section XIX.6]{Kassel}, define the algebra $S$ of formal series (in $\hbar$) of monomials in $A$ and $B$ and the $\bbC[[\hbar]]$-subalgebra $\bar{S}\subset S$ of those formal series whose monomials begin with $A$ and end with $B$, e.g. $1+\hbar^2AB$. Let $\pi:S\rightarrow\bar{S}$ denote the canonical projection then we can note that 
\begin{equation}\label{eq:lim pi W^c_epsilon}
\lim_{\varepsilon\to0}\pi\left(W^{c_\varepsilon}\right)=1+\sum _{\begin{smallmatrix}g\geq1\\p_1,\ldots,p_g\geq1\\q_1,\ldots,q_g\geq1\end{smallmatrix}}\hbar^{\sum_{i=1}^g(p_i+q_i)}(-1)^{\sum_{j=1}^gq_j}\zeta_{p_1+1,\{1\}^{q_1-1},\ldots,p_g+1,\{1\}^{q_g-1}}A^{p_1}B^{q_1}\cdots A^{p_g}B^{q_g}\,.
\end{equation}
Consider now the algebra $S[\alpha,\beta]$ of polynomials in two commuting variables $\alpha$ and $\beta$ with coefficients in $S$. We can always rewrite a monomial in $S[\alpha,\beta]$ as $\beta^sM\alpha^r$, where $M$ is a monomial in $S$. Define a $\bbC[[\hbar]]$-linear map $f':S[\alpha,\beta]\rightarrow S$ by
\begin{equation}\label{eq:def of f'}
f'(\beta^sM\alpha^r):=B^sMA^r\quad.
\end{equation}
From \eqref{eq:def of f'} we construct the $\bbC[[\hbar]]$-linear endomorphism $f:S\rightarrow S$ as
\begin{equation}
f\big(p(A,B)\big):=f'\big(p(A-\alpha,B-\beta)\big)
\end{equation}
where $p(A,B)\in S$. Notice that
\begin{equation}\label{eq:f kills MA and BM}
f(BM)=f(MA)=0\quad.
\end{equation}
\begin{lem}
The Drinfeld associator \eqref{eq:potentially finite Drinfeld} is finite and given by 
\begin{equation}\label{eq:Phi=f(pi W)}
\Phi(A,B)=f\left(\lim_{\varepsilon\to0}\pi\left(W^{c_\varepsilon}\right)\right)\quad.
\end{equation}
\end{lem}
\begin{proof}
Let us apply $f$ to the definition \eqref{eq:potentially finite Drinfeld}, 
\begin{equation}
f\big(\Phi(A,B)\big)=\lim_{\varepsilon\to0}f\left(e^{-\hbar\ln(\varepsilon)B}W^{c_\varepsilon}e^{\hbar\ln(\varepsilon)A}\right)=\lim_{\varepsilon\to0}f\left(W^{c_\varepsilon}\right)=f\left(\lim_{\varepsilon\to0}\pi\left(W^{c_\varepsilon}\right)\right)
\end{equation}
where the second equality comes from \eqref{eq:f kills MA and BM} and the third equality comes from the same fact together with \eqref{eq:c_epsilon compact PT}, \eqref{eq:convergent iterated integral} and \eqref{eq:lim pi W^c_epsilon}. We now wish to show that 
\begin{equation}
\Phi(A,B)=f\big(\Phi(A,B)\big):=f'\big(\Phi(A-\alpha,B-\beta)\big)\quad.
\end{equation}
In fact, it suffices to show that $\Phi(A-\alpha,B-\beta)=\Phi(A,B)$ and for this we consider the connection
\begin{equation}
\Gamma(A-\alpha,B-\beta):=\hbar\left(\frac{A-\alpha}{x}+\frac{B-\beta}{x-1}\right)dx=\Gamma(A,B)+\Gamma(-\alpha,-\beta)\quad.
\end{equation}
We already know that
\begin{equation}
^{\Gamma(A-\alpha,B-\beta)}W^{c_\varepsilon}=e^{\hbar\ln(\varepsilon)(B-\beta)}\Phi_\varepsilon(A-\alpha,B-\beta)e^{-\hbar\ln(\varepsilon)(A-\alpha)}
\end{equation}
and we can use the factorisation of parallel transport as in \cite[Item iv.) of Proposition 2]{BRW} to see that
\begin{align}
^{\Gamma(A-\alpha,B-\beta)}W^{c_\varepsilon}&={}^{\Gamma(-\alpha,-\beta)}W^{c_\varepsilon}\,{}^{\Gamma(A,B)}W^{c_\varepsilon}\nn\\
&=e^{-\hbar\ln(\varepsilon)\beta}e^{-\hbar\ln(1-\varepsilon)(\alpha-\beta)}e^{\hbar\ln(\varepsilon)\alpha}e^{\hbar\ln(\varepsilon)B}\Phi_\varepsilon(A,B)e^{-\hbar\ln(\varepsilon)A}
\end{align}
where we used Remark \ref{rem:at least diverging terms} in the second equality. Equating the above two expressions yields
\begin{equation}
\Phi_\varepsilon(A-\alpha,B-\beta)=e^{-\hbar\ln(1-\varepsilon)(\alpha-\beta)}\Phi_\varepsilon(A,B)\quad\implies\quad\Phi(A-\alpha,B-\beta)=\Phi(A,B)~.
\end{equation}
\end{proof}
We can explicitly compute the RHS of \eqref{eq:Phi=f(pi W)} as in \cite[Theorem A.9]{LeMu} by binomially expanding each term such as $(A-\alpha)^{p_1}$ and then applying \eqref{eq:def of f'}.
\begin{cor}\label{cor:Phi with MZV coeffs}
The Drinfeld associator \eqref{eq:Phi=f(pi W)} is explicitly given by
\begin{align}
&\Phi(A,B)=1+\sum_{r=2}^\infty\hbar^r\sum_{\begin{smallmatrix}\{g,p_1,q_1,\ldots,p_g,q_g\in\bbN^*\,|\,\sum_{i=1}^g(p_i+q_i)=r\}\end{smallmatrix}}(-1)^{\sum_{j=1}^gq_j}\,\zeta_{p_1+1,\{1\}^{q_1-1},\ldots,p_g+1,\{1\}^{q_g-1}}\nn\\
&\sum_{\begin{smallmatrix}0\leq s_1\leq p_1\,,\,0\leq t_1\leq q_1\\\vdots\\0\leq s_g\leq p_g\,,\,0\leq t_g\leq q_g\end{smallmatrix}}\left[\prod_{i=1}^g(-1)^{s_i+t_i}\binom{p_i}{s_i}\binom{q_i}{t_i}\right]B^{\sum_{n=1}^gt_n}A^{p_1-s_1}B^{q_1-t_1}\cdots A^{p_g-s_g}B^{q_g-t_g}A^{\sum_{m=1}^gs_m}
\end{align}
\end{cor}

\subsection{With coefficients as delta values}\label{subsec:BRW PT}
Let us recall the outline of the second method of constructing Drinfeld's Knizhnik-Zamolodchikov associator series. Going back to \eqref{eq:c_epsilon PT}, we factorise a continuous piecewise smooth reparametrisation of the affine path $c_\varepsilon(s):=\varepsilon+s(1-2\varepsilon)$ as the concatenation of two ``exponential half-paths". Explicitly, introduce:
\begin{equation}
c_{(1,\varepsilon)}(s):=\tfrac{1}{2}(2\varepsilon)^{1-s}\qquad,\qquad c_{(2,\varepsilon)}(s):=1-\tfrac{1}{2}(2\varepsilon)^s
\end{equation}
and 
\begin{equation}
\gamma(s):=\begin{cases}
    \frac{\tfrac{1}{2}(2\varepsilon)^{1-2s}-\varepsilon}{1-2\varepsilon}\quad\,\,,\qquad&0\leq s\leq\tfrac{1}{2}\\
    \frac{1-\tfrac{1}{2}(2\varepsilon)^{2s-1}-\varepsilon}{1-2\varepsilon}\,\,,\qquad&\tfrac{1}{2}\leq s\leq1
    \end{cases}
\end{equation}
then it is straightforward to check that $c_\varepsilon\circ\gamma=c_{(2,\varepsilon)}c_{(1,\varepsilon)}$. Furthermore, the \textbf{interval inversion} 
\begin{align}
\iota:[0,1]&\longrightarrow[0,1]\nn\\
s&\longmapsto1-s
\end{align} 
is such that
\begin{equation}\label{eq:nice facts about iota}
c_{(1,\varepsilon)}=\iota\circ c_{(2,\varepsilon)}\circ\iota\qquad,\qquad\iota^*\left(\Gamma(A,B)\right)=\Gamma(B,A)\quad.
\end{equation}
These facts (together with the standard functoriality properties of parallel transport \cite[Theorem 9]{BRW}) give us 
\begin{equation}
^{\Gamma(A,B)}W^{c_\varepsilon}={}^{\Gamma(A,B)}W^{c_{(2,\varepsilon)}}\left({}^{\Gamma(B,A)}W^{c_{(2,\varepsilon)}}\right)^{-1}\quad.
\end{equation}
In particular, \cite[Lemma 19]{BRW} tells us that 
\begin{subequations}
\begin{equation}
{}^{\Gamma(A,B)}W^{c_{(2,\varepsilon)}}=e^{\hbar\ln(\varepsilon)B}e^{\hbar\ln(2)B}\Xi^\varepsilon_{B,A}
\end{equation}
with
\begin{equation}
\Xi^\varepsilon_{B,A}:=1+\sum_{r=1}^\infty\hbar^r\sum_{\ell_1,\ldots,\ell_r=0}^\infty\hbar^{\sum_{i=1}^r\ell_i}\,\I^\varepsilon_{\ell_1\ldots\ell_r}\ad^{\ell_1}_B(A)\cdots\ad^{\ell_r}_B(A)
\end{equation}
and where
\begin{equation}\label{eq:epsilon iterated integral for associator}
\I^\varepsilon_{\ell_1\ldots\ell_r}:=\frac{1}{\ell_1!\cdots\ell_r!}\int_0^{-\ln(2\varepsilon)}\frac{\tau_1^{\ell_1}}{2e^{\tau_1}-1}d\tau_1\cdots\frac{\tau_r^{\ell_r}}{2e^{\tau_r}-1}d\tau_r\qquad.
\end{equation}
\end{subequations}
The proof of \cite[Lemma 19]{BRW} shows that \eqref{eq:epsilon iterated integral for associator} is indeed finite under the limit $\lim_{\varepsilon\to0}$. Therefore, \eqref{eq:potentially finite Drinfeld} tells us that\footnote{Notice as well that this form makes it obvious that $\Phi(A,B)=\Phi(B,A)^{-1}$.} 
\begin{subequations}
\begin{equation}\label{eq:associator as factor of PT}
\Phi(A,B)=e^{\hbar\ln(2)B}\Xi_{B,A}\big(\Xi_{A,B}\big)^{-1}e^{-\hbar\ln(2)A}
\end{equation}
with
\begin{equation}\label{eq:Xi_B,A:=}
\Xi_{B,A}:=1+\sum_{r=1}^\infty\hbar^r\sum_{\ell_1,\ldots,\ell_r=0}^\infty\hbar^{\sum_{i=1}^r\ell_i}\,\I_{\ell_1\ldots\ell_r}\ad^{\ell_1}_B(A)\cdots\ad^{\ell_r}_B(A)
\end{equation}
and where
\begin{equation}\label{eq:iterated integral for associator}
\I_{\ell_1\ldots\ell_r}:=\frac{1}{\ell_1!\cdots\ell_r!}\int_0^\infty\frac{\tau_1^{\ell_1}}{2e^{\tau_1}-1}d\tau_1\cdots\frac{\tau_r^{\ell_r}}{2e^{\tau_r}-1}d\tau_r\qquad.
\end{equation}
\end{subequations}
\begin{rem}\label{rem:iterated integral shuffle relations}
Notice that \eqref{eq:shuffle prod for iterated integrals} gives 
\begin{equation}\label{eq:shuffle for I}
\I_{\ell_1\ldots\ell_r}\I_{\ell_{r+1}\ldots\ell_{r+r'}}=\sum_{\sigma\in\mathrm{S}_{r,r'}}\I_{\ell_{\sigma(1)}\ldots\ell_{\sigma(r+r')}}\quad.
\end{equation}
We mention the following relevant instances of these shuffle product relations \eqref{eq:shuffle for I}:
\begin{subequations}
\begin{alignat}{3}
\I_i\I_j=\,&\I_{ij}+\I_{ji}\quad&&,\label{eq:I_iI_j}\\
\I_i\I_{jk}=\,&\I_{ijk}+\I_{jik}+\I_{jki}\quad&&,\label{eq:I_iI_jk}\\
\I_i\I_{jkl}=\,&\I_{ijkl}+\I_{jikl}+\I_{jkil}+\I_{jkli}\quad&&,\label{eq:I_iI_jkl}\\
\I_{ij}\I_{kl}=\,&\I_{ijkl}+\I_{ikjl}+\I_{iklj}+\I_{kijl}+\I_{kilj}+\I_{klij}\quad&&.\label{eq:I_ijI_jkl}
\end{alignat}
\end{subequations}
At third order, we will make use of
\begin{equation}\label{eq:third order shuffle relations}
\I_0\I_1=\I_{01}+\I_{10}\quad.
\end{equation}
At fourth order, we will use:
\begin{equation}\label{eq:fourth order shuffle relations}
\I_0\I_2=\I_{02}+\I_{20}\quad,\quad\I_0\I_{01}=2\I_{0^21}+\I_{010}\quad,\quad\I_1\I_{0^2}=\I_{10^2}+\I_{010}+\I_{0^21}\quad.
\end{equation}
At fifth order, we use:
\begin{alignat}{4}
\I_0\I_3=&\I_{03}+\I_{30}\quad&&,\quad\I_0\I_{02}=2\I_{0^22}+\I_{020}\quad&&&,\nn\\
\I_0\I_{0^21}=&3\I_{0^31}+\I_{0^210}\quad&&,\quad
\I_1\I_{0^3}=\I_{10^3}+\I_{010^2}+\I_{0^210}+\I_{0^31}\quad&&&,\nn\\
\I_2\I_{0^2}=&\I_{20^2}+\I_{020}+\I_{0^22}\quad&&,\quad\I_{0^2}\I_{01}=3\I_{0^31}+2\I_{0^210}+\I_{010^2}\quad&&&,\nn\\
\I_0\I_{1^2}=&2\I_{01^2}+\I_{101}+\I_{1^20}\quad&&&.
\end{alignat}
\end{rem}

\begin{propo}\label{propo:iterated integral is sum of deltas}
The iterated integrals $\I_{\ell_1\ldots\ell_r}$ of \eqref{eq:iterated integral for associator} are given by 
\begin{equation}\label{eq:I is sum of deltas}
\sum_{\begin{smallmatrix}0\leq m_1\leq\ell_1,\\
0\leq m_2\leq\ell_2+m_1,\\
\vdots\\
0\leq m_{r-1}\leq\ell_{r-1}+m_{r-2}\end{smallmatrix}}\binom{\ell_2+m_1}{\ell_2}\cdot\cdot\binom{\ell_r+m_{r-1}}{\ell_r}\delta_{\ell_r+m_{r-1}+1,\ell_{r-1}+m_{r-2}-m_{r-1}+1,\ldots,\ell_2+m_1-m_2+1,\ell_1-m_1+1}\,.
\end{equation}
\end{propo}
\begin{proof}
Let us make the substitution $\tau_i=-\ln x_i$ into \eqref{eq:iterated integral for associator},
\begin{align}
\I_{\ell_1\ldots\ell_r}=&\,\frac{(-1)^r}{\ell_1!\cdots\ell_r!}\int_1^0\frac{\left(\ln\frac{1}{x_1}\right)^{\ell_1}}{2-x_1}dx_1\frac{\left(\ln\frac{1}{x_2}\right)^{\ell_2}}{2-x_2}dx_2\cdots\frac{\left(\ln\frac{1}{x_{r-1}}\right)^{\ell_{r-1}}}{2-x_{r-1}}dx_{r-1}\frac{\left(\ln\frac{1}{x_r}\right)^{\ell_r}}{2-x_r}dx_r\nn\\
=&\,\frac{1}{\ell_1!\cdots\ell_r!}\int_0^1\frac{\left(\ln\frac{1}{x_r}\right)^{\ell_r}}{2-x_r}dx_r\frac{\left(\ln\frac{1}{x_{r-1}}\right)^{\ell_{r-1}}}{2-x_{r-1}}dx_{r-1}\cdots\frac{\left(\ln\frac{1}{x_2}\right)^{\ell_2}}{2-x_2}dx_2\frac{\left(\ln\frac{1}{x_1}\right)^{\ell_1}}{2-x_1}dx_1\quad,\label{eq:inverted iterated integral BRW}
\end{align}
where we used the inversion for iterated integrals \eqref{eq:inversion for iterated integrals} in the second equality. Let us binomially expand the last integrand as follows,
\begin{equation}
\left(\ln\tfrac{1}{x_1}\right)^{\ell_1}=\left(\ln\tfrac{x_2}{x_1}+\ln\tfrac{1}{x_2}\right)^{\ell_1}=\sum_{m_1=0}^{\ell_1}\binom{\ell_1}{m_1}\left(\ln\tfrac{1}{x_2}\right)^{m_1}\left(\ln\tfrac{x_2}{x_1}\right)^{\ell_1-m_1}\,.
\end{equation}
We perform another binomial expansion,
\begin{equation}
\left(\ln\tfrac{1}{x_2}\right)^{\ell_2+m_1}=\left(\ln\tfrac{x_3}{x_2}+\ln\tfrac{1}{x_3}\right)^{\ell_2+m_1}=\sum_{m_2=0}^{\ell_2+m_1}\binom{\ell_2+m_1}{m_2}\left(\ln\tfrac{1}{x_3}\right)^{m_2}\left(\ln\tfrac{x_3}{x_2}\right)^{\ell_2+m_1-m_2}\,.
\end{equation}
Iterating on this procedure of binomial expansions, we rewrite \eqref{eq:inverted iterated integral BRW} as
\begin{align}
&\frac{1}{\ell_1!\cdots\ell_r!}\sum_{m_1=0}^{\ell_1}\binom{\ell_1}{m_1}\sum_{m_2=0}^{\ell_2+m_1}\binom{\ell_2+m_1}{m_2}\cdots\sum_{m_{r-1}=0}^{\ell_{r-1}+m_{r-2}}\binom{\ell_{r-1}+m_{r-2}}{m_{r-1}}\label{eq:compare iterated BRW}\\
&\int_0^1\frac{\left(\ln\frac{1}{x_r}\right)^{\ell_r+m_{r-1}}}{2-x_r}dx_r\frac{\left(\ln\tfrac{x_r}{x_{r-1}}\right)^{\ell_{r-1}+m_{r-2}-m_{r-1}}}{2-x_{r-1}}dx_{r-1}\cdots\frac{\left(\ln\tfrac{x_3}{x_2}\right)^{\ell_2+m_1-m_2}}{2-x_2}dx_2\frac{\left(\ln\tfrac{x_2}{x_1}\right)^{\ell_1-m_1}}{2-x_1}dx_1\nn\,.
\end{align}
For $\{\ell'_i\in\bbN\,|\,1\leq i\leq r\}$, we substitute $z=\frac{1}{2}$ and $\{s_i=\ell'_{r+1-i}+1\,|\,1\leq i\leq r\}$ into \eqref{eq:single variable mult polylog as iterated int},
\begin{equation}\label{eq:compare SVoMP}
\delta_{\ell'_r+1,\ldots,\ell'_1+1}=\frac{1}{\ell'_1!\cdots\ell'_r!}\int_0^1\frac{\left(\ln\frac{1}{y_1}\right)^{\ell'_r}}{2-y_1}dy_1\frac{\left(\ln\frac{y_1}{y_2}\right)^{\ell'_{r-1}}}{2-y_2}dy_2\cdots\frac{\left(\ln\frac{y_{r-1}}{y_r}\right)^{\ell'_1}}{2-y_r}dy_r\quad.
\end{equation}
Finally, we compare \eqref{eq:compare iterated BRW} with \eqref{eq:compare SVoMP}.
\end{proof}
If $\ell_1=\ell_2=\ldots=\ell_r=0$ then we use \eqref{eq:I is sum of deltas} together with \eqref{eq:delta_(1^k)} to see that
\begin{equation}\label{eq:I_(0)^r}
\I_{0^r}=\delta_{\{1\}^r}=\frac{(\ln2)^r}{r!}\quad.
\end{equation}
It is also worth mentioning the special cases: where $r=1$,
\begin{equation}\label{eq:I_ell}
\I_\ell=\delta_{\ell+1}
\end{equation}
and where $r=2$,
\begin{equation}\label{eq:I_ell1,ell2}
\I_{\ell_1\ell_2}=\sum_{m=0}^{\ell_1}\binom{\ell_2+m}{\ell_2}\delta_{\ell_2+m+1,\ell_1-m+1}\quad.
\end{equation}
Using \eqref{eq:I_ell} and \eqref{eq:I_ell1,ell2}, the simplest shuffle relation \eqref{eq:third order shuffle relations} reads as
\begin{equation}\label{eq:delta_1,2}
\delta_2\ln2=2\delta_{2,1}+\delta_{1,2}\quad\implies\quad\delta_{1,2}=\frac{\zeta_2\ln2}{2}-\frac{(\ln2)^3}{6}-\frac{\zeta_3}{4}
\end{equation}
where the latter equality comes from \eqref{delta_2,1} and Euler's dilogarithm formula \eqref{eq:Li_2(1/2)}.
\begin{rem}
Choosing $i=j=1$ in \eqref{eq:I_iI_j} together with \eqref{eq:I_ell} and \eqref{eq:I_ell1,ell2} gives
\begin{equation}\label{eq:delta_2^2}
\delta_2^2=4\delta_{3,1}+2\delta_{2,2}\quad,
\end{equation}
while $(i=1,j=2)$ in \eqref{eq:I_iI_j} gives
\begin{equation}\label{eq:delta2delta3}
\delta_2\delta_3=6\delta_{4,1}+3\delta_{3,2}+\delta_{2,3}\quad.
\end{equation}
These shuffle relations hold more generally for single-variable multiple polylogarithms and not just delta values. This can be shown by first defining, for $\ell_1,\ldots,\ell_r\in\bbN^*$,
\begin{equation}
L^z_{\ell_1\ldots\ell_r}:=\frac{1}{\ell_1!\cdots\ell_r!}\int_0^\infty\frac{\tau_1^{\ell_1}}{\tfrac{1}{z}e^{\tau_1}-1}d\tau_1\cdots\frac{\tau_r^{\ell_r}}{\tfrac{1}{z}e^{\tau_r}-1}d\tau_r\quad.
\end{equation}
In particular, setting $z=1$ gives
\begin{subequations}\label{subeq:zeta_2^2 and zeta_2*zeta_3}
\begin{equation}\label{eq:zeta2^2}
\frac{\zeta_2^2}{2}=2\zeta_{3,1}+\zeta_{2,2}\quad.
\end{equation}
and
\begin{equation}\label{eq:zeta2zeta3}
\zeta_2\zeta_3=6\zeta_{4,1}+3\zeta_{3,2}+\zeta_{2,3}\quad.
\end{equation}
\end{subequations}
\end{rem}
Substituting the value for $\delta_{2,2}$ in \eqref{eq:delta_2,2} and Euler's dilogarithm formula \eqref{eq:Li_2(1/2)} into the shuffle relation \eqref{eq:delta_2^2} gives
\begin{equation}\label{eq:delta_3,1}
\delta_{3,1}=\frac{\zeta_4}{8}-\frac{\zeta_3\ln2}{8}+\frac{(\ln2)^4}{24}\quad.
\end{equation}
Substituting \eqref{eq:zeta_4,1} and the ``stuffle identity" \eqref{eq:stuffle identity} into the shuffle relation \eqref{eq:zeta2zeta3} gives 
\begin{equation}\label{eq:zeta_3,2}
\zeta_{3,2}=3\zeta_3\zeta_2-\tfrac{11}{2}\zeta_5\quad.
\end{equation}

\section{Fifth-order expansion of Drinfeld's KZ associator}\label{sec:Fifth-order expansion}
For sake of more compact expressions, we absorb the factor of $\hbar$ into $A$ and $B$ throughout.
\subsection{With coefficients as MZVs}\label{subsec:MZV expansion}
In this subsection we express the Drinfeld associator series of Corollary \ref{cor:Phi with MZV coeffs} modulo $\hbar^6$. Throughout this subsection we will repeatedly use the fact that, for $m\in\bbN$, 
\begin{equation}\label{eq:ad_A^m(B)}
\ad_A^m(B)=\sum_{s=0}^m(-1)^s\binom{m}{s}A^{m-s}BA^s\quad.
\end{equation}
Using \eqref{eq:ad_A^m(B)} we can see that, for $p,q\in\bbN^*$,
\begin{subequations}
\begin{align}
(-1)^q\sum_{s=0}^p\sum_{t=0}^q(-1)^{s+t}\binom{p}{s}\binom{q}{t}B^tA^{p-s}B^{q-t}A^s=&\,\sum_{s=0}^{p-1}(-1)^s\binom{p}{s}\ad^q_B(A^{p-s})A^s\label{eq:p<q}\\
=&\,\sum_{t=0}^{q-1}(-1)^{q+t}\binom{q}{t}B^t\ad^p_A(B^{q-t})\quad.\label{eq:q<p}
\end{align}
\end{subequations}
We will use \eqref{eq:p<q} when $p\leq q$ and \eqref{eq:q<p} when $q<p$. We can also see that, for $p_1,q_1,p_2,q_2\in\bbN^*$,
\begin{align}
(-1)^{q_1+q_2}&\sum_{\begin{smallmatrix}0\leq s_1\leq p_1\,,\,0\leq t_1\leq q_1\\
0\leq s_2\leq p_2\,,\,0\leq t_2\leq q_2\end{smallmatrix}}\left(\prod_{i=1}^2(-1)^{s_i+t_i}\binom{p_i}{s_i}\binom{q_i}{t_i}\right)B^{t_1+t_2}A^{p_1-s_1}B^{q_1-t_1}A^{p_2-s_2}B^{q_2-t_2}A^{s_1+s_2}\nn\\
&=(-1)^{q_2}\sum_{s_1=0}^{p_1-1}\sum_{t_2=0}^{q_2-1}(-1)^{s_1+t_2}\binom{p_1}{s_1}\binom{q_2}{t_2}B^{t_2}\ad^{q_1}_B(A^{p_1-s_1})\ad^{p_2}_A(B^{q_2-t_2})A^{s_1}\quad.
\end{align}
For $r\geq2$, $g\geq1$ and $\{p_1,q_1,\ldots,p_g,q_g\in\bbN^*\,|\,\sum_{i=1}^g(p_i+q_i)=r\}$:
\begin{itemize}
\item If $r=2$ then $g=p_1=q_1=1$. 
\item If $r=3$ then $g=1$ and $(p_1=2,q_1=1)$/$(p_1=1,q_1=2)$.
\item If $r=4$ then $g=1$ and $(p_1=3,q_1=1)$/$p_1=q_1=2$/$(p_1=1,q_1=3)$ or $g=2$ and $p_1=q_1=p_2=q_2=1$.
\item If $r=5$ then $g=1$ and $(p_1=4,q_1=1)$/$(p_1=3,q_1=2)$/$(p_1=2,q_1=3)$/$(p_1=1,q_1=4)$ or $g=2$ and $(p_1=2,q_1=1,p_2=1,q_2=1)/(p_1=1,q_1=2,p_2=1,q_2=1)/(p_1=1,q_1=1,p_2=2,q_2=1)/(p_1=1,q_1=1,p_2=1,q_2=2)$.
\end{itemize}
Using this solution set, we have the following expression for the Drinfeld associator series of Corollary \ref{cor:Phi with MZV coeffs} modulo $\hbar^6$,
\begin{align}
&1+\zeta_2\ad_B(A)-\zeta_3\ad^2_A(B)+\zeta_{2,1}\ad^2_B(A)-\zeta_4\ad^3_A(B)+\zeta_{3,1}\left[\ad^2_B(A^2)-2\ad^2_B(A)A\right]\label{eq:prelim 5th order MZV expansion}\\
&+\zeta_{2,1,1}\ad^3_B(A)-\zeta_{2,2}\ad_B(A)\ad_A(B)-\zeta_5\ad^4_A(B)+\zeta_{4,1}\left[\ad_A^3(B^2)-2B\ad_A^3(B)\right]\nn\\&+\zeta_{3,1,1}\left[\ad_B^3(A^2)-2\ad_B^3(A)A\right]+\zeta_{2,1,1,1}\ad^4_B(A)-\zeta_{3,2}\left[\ad_B(A^2)\ad_A(B)-2\ad_B(A)\ad_A(B)A\right]\nn\\&-\zeta_{2,1,2}\ad_B^2(A)\ad_A(B)-\zeta_{2,3}\ad_B(A)\ad_A^2(B)+\zeta_{2,2,1}\left[\ad_B(A)\ad_A(B^2)-2B\ad_B(A)\ad_A(B)\right]\quad.\nn
\end{align}
Repeated application of the Leibniz rule gives: 
\begin{subequations}
\begin{alignat}{3}
\ad^2_B(A^2)=&A\ad^2_B(A)+2\ad_B(A)\ad_B(A)+\ad^2_B(A)A\quad&&,\\
\ad^3_B(A^2)=&A\ad^3_B(A)+3\ad_B(A)\ad^2_B(A)+3\ad^2_B(A)\ad_B(A)+\ad^3_B(A)A\quad&&.
\end{alignat}
\end{subequations}
We denote $[B,A]=:X$ and use the MZV duality \eqref{eq:relevant MZV dualities} together with the shuffle relations \eqref{subeq:zeta_2^2 and zeta_2*zeta_3} to write the fifth-order expansion \eqref{eq:prelim 5th order MZV expansion} as
\begin{align}
&1+\zeta_2X+\zeta_3[A+B,X]+\zeta_4\left(\ad^2_A+\ad^2_B\right)(X)+\zeta_{3,1}\big[A,[B,X]\big]+\frac{\zeta_2^2}{2}X^2+\zeta_5\left(\ad^3_A+\ad^3_B\right)(X)\nn\\&+\zeta_{4,1}\left(\big[B,\ad^2_A(X)\big]+\big[A,\ad^2_B(X)\big]\right)+\zeta_2\zeta_3\left(X[A,X]+[B,X]X\right)+(3\zeta_{4,1}+\zeta_{3,2})\ad^2_X(A-B)\,.
\end{align}
Finally, we eliminate the terms of depth 2 by using \eqref{eq:zeta_3,1}, \eqref{eq:zeta_4,1} and \eqref{eq:zeta_3,2},
\begin{align}
&1+\zeta_2X+\zeta_3[A+B,X]+\zeta_4\Big(\big(\ad^2_A+\ad^2_B\big)(X)+\tfrac{1}{4}\big[A,[B,X]\big]\Big)+\frac{\zeta_2^2}{2}X^2+\zeta_3\zeta_2\left(X[A,X]+[B,X]X\right)\nn\\&+\zeta_5\Big(\big(\ad^3_A+\ad^3_B\big)(X)+\tfrac{1}{2}\ad^2_X(A-B)\Big)+(2\zeta_5-\zeta_3\zeta_2)\left(\big[B,\ad^2_A(X)\big]+\big[A,\ad^2_B(X)\big]\right)\quad.\label{eq:fifth order MZV expansion}
\end{align}

\subsection{With coefficients as delta values}\label{subsec:DV expansion}
This subsection expresses the Drinfeld series \eqref{eq:associator as factor of PT} modulo $\hbar^6$ and denotes $\ln2=:c$ and $c(A+B)=:C$. With this notation, the Drinfeld series \eqref{eq:associator as factor of PT} is written as
\begin{equation}
\Phi(A,B)=e^{cB}\Xi_{B,A}\big(e^{cA}\Xi_{A,B}\big)^{-1}
\end{equation}
thus we denote
\begin{equation}\label{eq:psi^BA}
\psi^{BA}:=e^{cB}\Xi_{B,A}=1+C+\sum_{k=2}^\infty\psi^{BA}_k
\end{equation}
where $\psi^{BA}_k$ is the term of order $\hbar^k$. Consider
\begin{equation}
(\psi^{AB})^{-1}=1+\sum_{l=1}^\infty(-1)^l\left[C+\sum_{k=2}^\infty\psi^{AB}_k\right]^l\quad,
\end{equation}
This gives the following expression for $\Phi=\psi^{BA}(\psi^{AB})^{-1}$,
\begin{align}
&1+C+\sum_{k=2}^\infty\psi^{BA}_k+\sum_{l=1}^\infty(-1)^l\left[C+\sum_{k=2}^\infty\psi^{AB}_k\right]^l+\left[C+\sum_{k=2}^\infty\psi^{BA}_k\right]\sum_{l=1}^\infty(-1)^l\left[C+\sum_{k=2}^\infty\psi^{AB}_k\right]^l\nn\\
&=1+\sum_{k=2}^\infty(\psi^{BA}_k-\psi^{AB}_k)+\sum_{l=2}^\infty(-1)^l\left[C+\sum_{k=2}^\infty\psi^{AB}_k\right]^l+\left[C+\sum_{k=2}^\infty\psi^{BA}_k\right]\sum_{l=1}^\infty(-1)^l\left[C+\sum_{k=2}^\infty\psi^{AB}_k\right]^l\nn\\
&=1+\sum_{k=2}^\infty(\psi^{BA}_k-\psi^{AB}_k)\sum_{l=0}^\infty(-1)^l\left[C+\sum_{k=2}^\infty\psi^{AB}_k\right]^l\quad,
\end{align}
where we used $\sum_{l=2}^\infty(-1)^l\left[C+\sum_{k=2}^\infty\psi^{AB}_k\right]^l=-\left[C+\sum_{k=2}^\infty\psi^{AB}_k\right]\sum_{l=1}^\infty(-1)^l\left[C+\sum_{k=2}^\infty\psi^{AB}_k\right]^l$ for the third equality. Let us express the fifth-order truncation,
\begin{align}
&1+\psi^{BA}_2-\psi^{AB}_2+\psi^{BA}_3-\psi^{AB}_3-(\psi^{BA}_2-\psi^{AB}_2)C+\psi^{BA}_4-\psi^{AB}_4-(\psi^{BA}_3-\psi^{AB}_3)C\nn\\
&+(\psi^{BA}_2-\psi^{AB}_2)(C^2-\psi^{AB}_2)+\psi^{BA}_5-\psi^{AB}_5-(\psi^{BA}_4-\psi^{AB}_4)C+(\psi^{BA}_3-\psi^{AB}_3)(C^2-\psi^{AB}_2)\nn\\
&-(\psi^{BA}_2-\psi^{AB}_2)(\psi^{AB}_3-C\psi^{AB}_2-\psi^{AB}_2C+C^3)\quad.
\end{align}
Defining the \textbf{antisymmetrisation} $\omega_i:=\psi^{BA}_i-\psi^{AB}_i$, we rewrite the fifth-order truncation as 
\begin{equation}\label{eq:fifth order truncation BRW}
\Phi=1+\omega_2+(\omega_3-\Phi_2C)+(\omega_4-\Phi_3C-\Phi_2\psi^{AB}_2)+(\omega_5-\Phi_4C-\Phi_3\psi^{AB}_2-\Phi_2\psi^{AB}_3)~.
\end{equation}
Using \eqref{eq:psi^BA}, \eqref{eq:Xi_B,A:=} and \eqref{eq:I_(0)^r} gives
\begin{equation}
\psi^{BA}_2=\frac{c^2}{2}(A^2+2BA+B^2)+\I_1X=\frac{C^2}{2}+\left(\frac{c^2}{2}+\I_1\right)X
\end{equation}
thus
\begin{equation}\label{eq:2nd order BRW}
\Phi_2=\omega_2=(c^2+2\I_1)X\quad.
\end{equation}
In particular, comparing \eqref{eq:2nd order BRW} with the second-order term of \eqref{eq:fifth order MZV expansion} yields
\begin{equation}\label{delta2 zeta2 relation}
c^2+2\I_1=\zeta_2\quad.
\end{equation}
We also have
\begin{align}
\psi^{BA}_3=\,&\frac{c^3}{6}(A^3+3BA^2+3B^2A+B^3)+c\I_1BX+\I_2[B,X]+\I_{10}XA+\I_{01}AX\nn\\
=\,&\frac{C^3}{6}+\frac{\zeta_2}{2}XC+\left(\I_2+c\I_1+\frac{c^3}{3}\right)[B,X]+\left(\I_{01}+\frac{c^3}{6}\right)[A,X]
\end{align}
where the second equality used \eqref{eq:third order shuffle relations} and \eqref{delta2 zeta2 relation}. The third-order antisymmetrisation is thus
\begin{align}
\omega_3=\zeta_2XC+\left(\I_1+\frac{c^2}{2}+\frac{\I_2+\I_{01}}{c}\right)[C,X]
\end{align}
and the third-order terms become 
\begin{equation}\label{eq:3rd order BRW}
\Phi_3=\omega_3-\zeta_2XC=\left(\I_1+\frac{c^2}{2}+\frac{\I_2+\I_{01}}{c}\right)[C,X]\quad.
\end{equation}
Comparing \eqref{eq:3rd order BRW} with the third-order term of \eqref{eq:fifth order MZV expansion} gives
\begin{equation}\label{eq:third order relation compare}
c\I_1+\frac{c^3}{2}+\I_2+\I_{01}=\,\zeta_3\quad.
\end{equation}
The fourth-order terms are thus given by
\begin{equation}\label{eq:prelim Phi_4}
\Phi_4=\omega_4-\Phi_3C-\Phi_2\psi^{AB}_2=\omega_4-\frac{\zeta_3}{c}[C,X]C-\frac{\zeta_2}{2}XC^2+\frac{\zeta_2^2}{2}X^2\quad.
\end{equation}
We have the following expression for $\psi_4^{BA}$,
\begin{align}
&\,\frac{c^4}{24}(A^4+4BA^3+6B^2A^2+4B^3A+B^4)+\I_{10^2}XA^2+\I_{010}AXA+\I_{0^21}A^2X+\I_{1^2}X^2\nn\\&+\I_{20}[B,X]A+\I_{02}A[B,X]+\I_3\ad^2_B(X)+c\I_{10}BXA+c\I_{01}BAX+c\I_2B[B,X]+\frac{c^2}{2}\I_1B^2X
\end{align}
thus we use \eqref{eq:third order shuffle relations}, \eqref{eq:fourth order shuffle relations}, \eqref{delta2 zeta2 relation} and \eqref{eq:third order relation compare} to express the fourth-order antisymmetrisation as
\begin{align}
\omega_4=\,&\frac{\zeta_2}{2}XC^2+\left(2(\I_{02}+c\I_{01})+\frac{c^4}{4}\right)\big[A,[B,X]\big]\nn\\&+\left(\frac{c^4}{6}+\frac{c^2}{2}\I_1+\I_3+c\I_2+\I_{0^21}\right)\left(\ad^2_B+\ad^2_A\right)(X)+\frac{\zeta_3}{c}[C,X]C\quad.
\end{align}
We thus have the following expression for the fourth-order terms \eqref{eq:prelim Phi_4},
\begin{align}\label{eq:Phi_4 expression}
\left(\frac{c^4}{6}+\frac{c^2}{2}\I_1+\I_3+c\I_2+\I_{0^21}\right)\left(\ad^2_B+\ad^2_A\right)(X)+\frac{\zeta_2^2}{2}X^2+\left(2(\I_{02}+c\I_{01})+\frac{c^4}{4}\right)\big[A,[B,X]\big]
\end{align}
and comparing these fourth-order terms \eqref{eq:Phi_4 expression} with the fourth-order terms of \eqref{eq:fifth order MZV expansion} yields:
\begin{subequations}\label{sube:fourth order compare}
\begin{alignat}{3}
\frac{c^4}{6}+\frac{c^2}{2}\I_1+\I_3+c\I_2+\I_{0^21}=&\zeta_4\quad&&,\label{eq:n=2 Holder}\\
2(\I_{02}+c\I_{01})+\frac{c^4}{4}=&\tfrac{1}{4}\zeta_4\quad&&.\label{eq:not actually new relation}
\end{alignat}
\end{subequations}
Substituting these relations into the fifth-order terms of \eqref{eq:fifth order truncation BRW}, $\Phi_5$ is expressed as
\begin{align}
&\omega_5-\tfrac{1}{4}\zeta_4\big[A,[B,X]\big]C-\zeta_4\left(\ad^2_B+\ad^2_A\right)(X)C-\frac{\zeta_3}{2c}[C,X]C^2-\frac{\zeta_2}{6}XC^3\\&+\zeta_2\zeta_3\left(X[A,X]+[B,X]X\right)+\zeta_2\left(\frac{1}{2}\zeta_3-\frac{c^3}{6}-\I_{01}\right)[A-B,X]X+\zeta_2\left(\frac{c^3}{6}+\I_{01}\right)\ad^2_X(A-B)\,.\nn
\end{align}
Finally, $\psi^{BA}_5$ is given by
\begin{align}
&\frac{c^5}{120}(A^5+5BA^4+10B^2A^3+10B^3A^2+5B^4A+B^5)+\I_{10^3}XA^3+\I_{010^2}AXA^2+\I_{0^210}A^2XA\nn\\&+\I_{0^31}A^3X+\I_{1^20}X^2A+\I_{101}XAX+\I_{01^2}AX^2+\I_{20^2}[B,X]A^2+\I_{020}A[B,X]A+\I_{0^22}A^2[B,X]\nn\\&+\I_{30}\ad^2_B(X)A+\I_{03}A\ad^2_B(X)+\I_{21}[B,X]X+\I_{12}X[B,X]+\I_4\ad^3_B(X)+c\I_{10^2}BXA^2\nn\\&+c\I_{010}BAXA+c\I_{0^21}BA^2X+c\I_{1^2}BX^2+c\I_{20}B[B,X]A+c\I_{02}BA[B,X]+c\I_3B\ad^2_B(X)\nn\\&+\frac{c^2}{2}\I_{10}B^2XA+\frac{c^2}{2}\I_{01}B^2AX+\frac{c^2}{2}\I_2B^2[B,X]+\frac{c^3}{6}\I_1B^3X\quad.
\end{align}
After computing the fifth-order antisymmetrisation by using the shuffle relations of Remark \ref{rem:iterated integral shuffle relations} together with \eqref{delta2 zeta2 relation}, \eqref{eq:third order relation compare} and \eqref{sube:fourth order compare}, the fifth-order terms become
\begin{align}
&\left(\I_{03}+c\I_{02}+\frac{c^2}{2}\I_{01}+\I_{0^22}+c\I_{0^21}+\frac{c^5}{12}\right)\left([A,\ad^2_B(X)]+[B,\ad^2_A(X)]\right)\nn\\&+\left(\I_4+c\I_3+\frac{c^2}{2}\I_2+\frac{c^3}{6}\I_1+\I_{0^31}+\frac{c^5}{24}\right)\left(\ad^3_A+\ad^3_B\right)(X)+\zeta_2\zeta_3\left(X[A,X]+[B,X]X\right)\nn\\&+\left(\I_{12}+\frac{c}{2}\I_1^2+c\I_{02}+\left[\frac{c^2}{2}+\zeta_2\right]\I_{01}+\I_{0^22}-\I_{101}-\I_{01^2}+\frac{c^3}{4}\zeta_2\right)\ad^2_X(A-B)\quad.
\end{align}
Comparing these fifth-order terms with the ones in \eqref{eq:fifth order MZV expansion} yields:
\begin{subequations}
\begin{alignat}{3}
\I_{03}+c\I_{02}+\frac{c^2}{2}\I_{01}+\I_{0^22}+c\I_{0^21}+\frac{c^5}{12}=&\,2\zeta_5-\zeta_3\zeta_2\quad&&,\label{eq:new relation 1}\\
\I_4+c\I_3+\frac{c^2}{2}\I_2+\frac{c^3}{6}\I_1+\I_{0^31}+\frac{c^5}{24}=&\,\zeta_5\quad&&,\\
\I_{12}+\frac{c}{2}\I_1^2+c\I_{02}+\left[\frac{c^2}{2}+\zeta_2\right]\I_{01}+\I_{0^22}-\I_{101}-\I_{01^2}+\frac{c^3}{4}\zeta_2=&\,\tfrac{1}{2}\zeta_5\quad&&.\label{eq:new relation 2}
\end{alignat}
\end{subequations}

\subsection{Relations from equating the fifth-order expansions}\label{subsec:relations}
Using \eqref{eq:I_ell}, the second-order relation \eqref{delta2 zeta2 relation} reads
\begin{equation}\label{eq:Holder n=0}
c^2+2\delta_2=\zeta_2\quad,
\end{equation}
i.e. \eqref{eq:Holder convolution} with $n=0$. Using \eqref{eq:I_ell1,ell2}, the third-order relation \eqref{eq:third order relation compare} reads
\begin{equation}\label{eq:Holder n=1}
\delta_{2,1}+\frac{c^3}{2}+c\delta_2+\delta_3=\zeta_3\quad,
\end{equation}
i.e. \eqref{eq:Holder convolution} with $n=1$. Using \eqref{eq:I is sum of deltas}, the fourth-order relations \eqref{sube:fourth order compare} read:
\begin{subequations}
\begin{alignat}{3}
\delta_{2,1,1}+\frac{c^4}{6}+\frac{c^2}{2}\delta_2+c\delta_3+\delta_4=&\,\zeta_4\quad&&,\label{eq:nonmysterious relation}\\
2c\delta_{2,1}+2\delta_{3,1}+\frac{c^4}{4}=&\,\tfrac{1}{4}\zeta_4\quad&&.\label{eq:mysterious relation}
\end{alignat}
\end{subequations}
The relation \eqref{eq:nonmysterious relation} is \eqref{eq:Holder convolution} with $n=2$ whereas \eqref{eq:mysterious relation} derives from \eqref{delta_2,1} and \eqref{eq:delta_3,1}.
\begin{rem}\label{rem:genuinely new}
Let us now focus on the fifth-order relations \eqref{eq:new relation 1} and \eqref{eq:new relation 2}:
\begin{subequations}\label{sube:new relations updated}
\begin{alignat}{2}
2\zeta_5-\zeta_2\zeta_3=&\,\delta_{4,1}+\frac{c\zeta_4}{8}-\frac{c^2\zeta_3}{16}+\frac{c^5}{24}+\delta_{3,1,1}+c\delta_{2,1,1}&&,\\
\tfrac{1}{2}\zeta_5-\delta_{3,2}-3\delta_{4,1}=&\,\frac{c\zeta_2^2}{8}+\frac{c\zeta_4}{8}+\frac{c^5}{12}-\frac{c^2\zeta_3}{16}+\frac{\zeta_2\zeta_3}{8}-\frac{c^3\zeta_2}{6}-3\delta_{3,1,1}-2\delta_{2,2,1}-\delta_{2,1,2}~&&.
\end{alignat}
\end{subequations}
Let us consider the following four shuffle relations coming from \eqref{eq:I_iI_jk}:
\begin{subequations}
\begin{alignat}{4}
c\delta_{3,1}&=3\delta_{3,1,1}+\delta_{2,2,1}+\delta_{1,3,1}\qquad&&,\\
\delta_2\delta_{2,1}&=6\delta_{3,1,1}+3\delta_{2,2,1}+\delta_{2,1,2}\qquad&&,\\
c\delta_{2,2}-\delta_{1,2,2}&=2\delta_{2,2,1}+\delta_{2,1,2}\qquad&&,\\
\delta_2\delta_{1,2}-2\delta_{1,2,2}&=2\delta_{2,2,1}+4\delta_{1,3,1}+2\delta_{2,1,2}\qquad&&.
\end{alignat}
\end{subequations}
We invert this system of equalities to get:
\begin{subequations}\label{sube:inverted system of equalities}
\begin{alignat}{3}
\delta_{3,1,1}&=-\tfrac{1}{6}c\delta_{3,1}+\tfrac{1}{4}\delta_2\delta_{2,1}-\tfrac{1}{3}c\delta_{2,2}+\tfrac{1}{24}\delta_2\delta_{1,2}+\tfrac{1}{4}\delta_{1,2,2}\qquad&&,\\
\delta_{2,2,1}&=c\delta_{3,1}-\tfrac{1}{2}\delta_2\delta_{2,1}+c\delta_{2,2}-\tfrac{1}{4}\delta_2\delta_{1,2}-\tfrac{1}{2}\delta_{1,2,2}\qquad&&,\\
\delta_{2,1,2}&=-2c\delta_{3,1}+\delta_2\delta_{2,1}-c\delta_{2,2}+\tfrac{1}{2}\delta_2\delta_{1,2}\qquad&&.
\end{alignat}
\end{subequations}
Using the first equality of \eqref{eq:delta_1,2}, we can rewrite \eqref{sube:inverted system of equalities} as:
\begin{subequations}
\begin{alignat}{3}
\delta_{3,1,1}&=\frac{1}{24}\delta_2\left(4\delta_{2,1}+c\delta_2\right)-\frac{c}{6}\left(\delta_{3,1}+2\delta_{2,2}\right)+\tfrac{1}{4}\delta_{1,2,2}\qquad&&,\label{eq:delta_3,1,1}\\
\delta_{2,2,1}&=\frac{1}{4}c\left(4\delta_{3,1}+4\delta_{2,2}-\delta_2^2\right)-\tfrac{1}{2}\delta_{1,2,2}\qquad&&,\\
\delta_{2,1,2}&=\frac{c}{2}\left(\delta_2^2-4\delta_{3,1}-2\delta_{2,2}\right)\qquad&&,
\end{alignat}
\end{subequations}
thus
\begin{equation}\label{eq:actual relevant linear combo}
3\delta_{3,1,1}+2\delta_{2,2,1}+\delta_{2,1,2}=\frac{1}{8}\delta_2\left(4\delta_{2,1}+c\delta_2\right)-\frac{1}{4}\delta_{1,2,2}-\frac{c}{2}\delta_{3,1}\quad.
\end{equation}
Substituting \eqref{eq:Li_2(1/2)}, \eqref{delta_2,1}, \eqref{eq:delta_3,1} and \eqref{subeq:delta_2,2 and delta_1,2,2} into \eqref{eq:delta_3,1,1} and \eqref{eq:actual relevant linear combo} gives:
\begin{subequations}\label{subeq:penultimate pair of relations}
\begin{alignat}{2}
\delta_{3,1,1}=&\,\frac{3\zeta_5}{64}-\frac{\zeta_2\zeta_3}{48}-\frac{c^2\zeta_3}{24}+\frac{c^5}{180}+\frac{c^3\zeta_2}{36}\quad&&,\\
3\delta_{3,1,1}+2\delta_{2,2,1}+\delta_{2,1,2}=&\,\frac{\zeta_2\zeta_3}{16}-\frac{c^3\zeta_2}{12}-\frac{3\zeta_5}{64}+\frac{c^5}{20}\quad&&.
\end{alignat}
\end{subequations}
Substituting \eqref{eq:nonmysterious relation} and \eqref{subeq:penultimate pair of relations} into \eqref{sube:new relations updated} gives:
\begin{subequations}\label{sube:revealed importance}
\begin{alignat}{2}
\delta_{4,1}=&\,c\delta_4+\frac{125\zeta_5}{64}+\frac{47c^2\zeta_3}{48}-\frac{47\zeta_2\zeta_3}{48}-\frac{9c\zeta_4}{8}-\frac{5c^3\zeta_2}{18}+\frac{13c^5}{360}\quad&&,\\
\delta_{3,2}+3\delta_{4,1}=&\,\frac{29\zeta_5}{64}-\frac{c\zeta_2^2}{8}-\frac{c\zeta_4}{8}-\frac{c^5}{30}+\frac{c^2\zeta_3}{16}-\frac{\zeta_2\zeta_3}{16}+\frac{c^3\zeta_2}{12}\quad&&.
\end{alignat}
\end{subequations}
We also have the following shuffle relation coming from \eqref{eq:I_iI_j},
\begin{equation}\label{eq:last shuffle relation}
c\delta_4=\,\delta_{1,4}+\delta_{2,3}+\delta_{3,2}+2\delta_{4,1}\quad,
\end{equation}
together with the shuffle relation \eqref{eq:delta2delta3}.
The point is that the two relations \eqref{sube:revealed importance} together with the shuffle relations \eqref{eq:last shuffle relation} and \eqref{eq:delta2delta3} allow us to eliminate the depth 2 terms ($\delta_{1,4},\,\delta_{2,3},\,\delta_{3,2}$ and $\delta_{4,1}$) in favour of products of depth 1 terms.
\end{rem}

\section*{Acknowledgments}
I would like to thank the anonymous reviewer for a plethora of useful comments that greatly improved the presentation and structure of this paper. Thanks also go to Alexander Schenkel and Robert Laugwitz for helpful discussions in understanding the material of \cite[Subsection 2.2]{BRW}. This work was supported by an EPSRC doctoral studentship (ref: 2747173).

\end{document}